\documentclass[12pt]{article}
\usepackage{amsmath,amsfonts}
\date{}

\begin{document}

\baselineskip 20pt
\title{On a problem of von Renteln }

\author{Arthur A.~Danielyan}

\maketitle

\begin{abstract}

\noindent We discuss the solution of a harmonic analysis problem proposed by M. von Renteln in 1980.

\end{abstract}

In 1980 Michael von Renteln has formulated the following Problem 5.62 and the remarks, which
have been published in the well-known collection (`Hayman's List')  of research problems
 in complex analysis
by D. Campbell, J. Clunie, and W. Hayman \cite{CCH} (p. 554).

\vspace{0.25 cm}
 
{\bf Problem 5.62.} Let $D= \{z: |z|<1 \} , \  T=\{z : |z|=1 \}$ and   $u$ be a continuous real-valued function on $T$. 
Give a necessary and sufficient condition on $u$ such that $u$ is the real part of a function 
$f$ in the disc algebra $A(\overline{D})$.

\vspace{0.15 cm}

{\it Remarks}. 

(1) A solution would have applications in the algebraic ideal theory.

(2) An answer to the analogous problem for $L^{p}(T)$,  $H^{p}(D)$ is the Burkholder-Gundy-Silverstein 
Theorem (see Peterson \cite{Pe}, p. 13).

\vspace{0.25 cm}

Note that in the recent account \cite{HL}  
of the mentioned collection of problems the authors indicate that no progress on Problem 5.62 has been reported (see \cite{HL}, p. 125).

The present paper clarifies that in fact a solution of Problem 5.62 is formulated in an exercise 
of Zygmund's book \cite{Zyg}. More precisely, the Exercise 5, part (a), on p. 180 in  \cite{Zyg} formulates a solution (theorem) with a reference to a 
report of M. Zamansky \cite{Za}. 
We use the results of the classical harmonic analysis from \cite{Zyg} to present a proof of Zamansky's theorem.
This gives the complete
solution of Problem 5.62,
which was assumed to be open for 42 years. Publishing this solution may also be helpful for 
possible applications in particular in algebraic ideal theory as stated in above Remark 1.

 As in \cite{Zyg}, 
we use the standard notations $ \tilde{f}(r,x)$ and $\tilde{f}(x)$
 for the conjugate functions in the unit disc and on the (interval $[0, 2\pi]$ of) real line, respectively. 
 We need the following classical theorem on the
boundary behavior of the conjugate harmonic function (see Theorem (7.20) in \cite{Zyg}, p. 103).

\vspace{0.25 cm}
 
{\bf Theorem A.}  {\it If $f$ is integrable and $F$ the indefinite integral of $f$, then

$$ \tilde{f}(r,x) -   \left ( -\frac{1}{\pi} \int_{1-r}^{\pi}[f(x+t) - f(x-t)] \frac{1}{2} \cot \frac{1}{2} t \ dt   \right ) \rightarrow 0 \ \ \ (r \rightarrow 1) $$   
at every point where $F$ is smooth, in particular where $f$ is continuous. If $f$ is everywhere continuous, the convergence is uniform.}

\vspace{0.25 cm}

By the definition (see the top of the same p. 103 in \cite{Zyg}), the function $F$ is smooth at the point $x$ if $F(x+t) + F(x-t) - 2F(x) = o(t).$
Since obviously the existence of the finite derivative of $F$ at $x$ implies the smoothness of $F$ at $x$,
the limit relation in Theorem A is true a. e.

The following theorem, which solves Problem 5.62, we formulate as in Exercise 5, part (a), on p. 180 in \cite{Zyg}.

\vspace{0.25 cm}

{\bf Zamansky's Theorem.}  {\it Let $f(x)$ be continuous and periodic. A necessary and sufficient condition for $\tilde{f}(x)$
to be continuous is that 

$$ \tilde{f}(x,h) = -\frac{1}{\pi} \int_{h}^{\pi}[f(x+t) - f(x-t)] \frac{1}{2} \cot \frac{1}{2} t \ dt    $$   
converges uniformly as $h \rightarrow +0.$}

\vspace{0.25 cm}

Before proving Zamansky's theorem, note that in the formulation of this theorem
we presented $ \tilde{f}(x,h)$ in slightly more explicit form than in \cite{Zyg}, p. 180. To see that the above form
is the same as that in \cite{Zyg}, we refer the reader to the simple notation adopted on p. 50 in \cite{Zyg}.

\vspace{0.25 cm}

{\bf Proof.} {(\it Sufficiency.)} Let  $\tilde{f}(x,h)$ be uniformly convergent as $h \rightarrow +0.$ Since $f$ is continuous, $f \in L$ 
 and the function  
$$ \tilde{f}(x) = -\frac{1}{\pi} \int_{0}^{\pi}[f(x+t) - f(x-t)] \frac{1}{2} \cot \frac{1}{2} t \ dt =  -\frac{1}{\pi} \lim_{h \rightarrow 0}  \int_{h}^{\pi} $$    
exists for almost all $x$ (see Theorem (3.1), p. 131, \cite{Zyg}).  By the definition of $ \tilde{f}(x,h)$,
the previous equality can be written as 
$\tilde{f}(x) = \lim_{h \rightarrow 0}\tilde{f}(x,h)$ for almost all $x$.  Since $\tilde{f}(x,h)$ converges uniformly, $\tilde{f}(x)$ is continuous.
   
 {\it (Necessity.)} Let $\tilde{f}(x)$ be continuous. We write the left part of the limit relation in Theorem A as follows:
$$ [ \tilde{f}(r,x) - \tilde{f}(x)] + \left [ \tilde{f}(x) -  \left ( -\frac{1}{\pi} \int_{1-r}^{\pi}[f(x+t) - f(x-t)] \frac{1}{2} \cot \frac{1}{2} t \ dt   \right ) \right] \rightarrow 0 \ \ \ (r \rightarrow 1)$$  
Since $f$ is continuous, by Theorem A the convergence in the last limit relation is uniform. 
Since 
$\tilde{f}(x)$ is continuous, we have that $ [ \tilde{f}(r,x) - \tilde{f}(x)] \rightarrow 0 \ \ \  (r \rightarrow 1)$ uniformly. 
Thus $$\left [ \tilde{f}(x) -  \left ( -\frac{1}{\pi} \int_{1-r}^{\pi}[f(x+t) - f(x-t)] \frac{1}{2} \cot \frac{1}{2} t \ dt   \right ) \right] \rightarrow 0 \ \ \ (r \rightarrow 1)$$ uniformly too. 
Because $\tilde{f}(x)$ does not depend on $r$, the convergence of the term 
$$ \left ( -\frac{1}{\pi} \int_{1-r}^{\pi}[f(x+t) - f(x-t)] \frac{1}{2} \cot \frac{1}{2} t \ dt   \right ) \ \ \ (r \rightarrow 1)$$
is uniform, which is the same as
$ \tilde{f}(x,h) = -\frac{1}{\pi} \int_{h}^{\pi}[f(x+t) - f(x-t)] \frac{1}{2} \cot \frac{1}{2} t \ dt $ converges uniformly as $h \rightarrow +0.$ This completes the proof.

\begin{minipage}[t]{6.5cm}
Arthur A. Danielyan\\
Department of Mathematics and Statistics\\
University of South Florida\\
Tampa, Florida 33620\\
USA\\
{\small e-mail: adaniely@usf.edu}
\end{minipage}

\end{document}